\documentclass{amsart}

\makeatletter
\let\@wraptoccontribs\wraptoccontribs
\makeatother

\usepackage{amsmath}
\usepackage{amsxtra}
\usepackage{amscd}
\usepackage{amsthm}
\usepackage{amsfonts}
\usepackage{amssymb}
\usepackage{mathtools}
\usepackage{eucal}
\usepackage[all]{xy}
\usepackage{graphicx}
\usepackage{tikz-cd}
\usepackage{mathrsfs}
\usepackage{euler}
\usepackage{color}
\usepackage{longtable}
\usepackage{caption}
\usepackage{enumerate}
\usepackage{faktor}
\usepackage{bbm}
\usepackage{stmaryrd}
\usepackage{rotating}
\usepackage[T1]{fontenc} % improved font encoding
\usepackage[utf8]{inputenc} % for better handling of non-ASCII characters
\usepackage{tikz}
\usepackage{tikz-cd}
\usetikzlibrary{matrix,arrows,decorations.pathmorphing}
\usepackage{hyperref}
\usepackage{float}
\usepackage{cleveref}

\usepackage[color=blue!10]{todonotes}
\theoremstyle{plain}
\newtheorem{theorem}{Theorem}
\newtheorem{lemma}[theorem]{Lemma}
\newtheorem{corollary}[theorem]{Corollary}
\theoremstyle{definition}

\theoremstyle{remark}
\newtheorem{remark}[theorem]{Remark}
\newtheorem{ex}[theorem]{Example}
\numberwithin{theorem}{section}
\numberwithin{equation}{section}

\DeclareMathAlphabet{\mathcal}{OMS}{cmsy}{m}{n}

\newcommand{\C}{\mathbb{C}}

\newcommand{\F}{\mathbb{F}}
\newcommand{\G}{\mathbb{G}}

\newcommand{\Q}{\mathbb{Q}}

\newcommand{\Z}{\mathbb{Z}}

\newcommand{\ttmat}[4]{\left( \begin{array}{cc}
#1 & #2 \\
#3 & #4
\end{array}
\right)}

\newcommand{\Ind}{\mathrm{Ind}}

\newcommand{\ov}[1]{{\overline{#1}}}
\newcommand{\onto}{\twoheadrightarrow}
\newcommand{\into}{\hookrightarrow}
\newcommand{\cO}{\mathcal{O}}
\DeclareMathOperator{\GL}{GL}
\newcommand{\cC}{\mathcal{C}}

\newcommand{\ann}{\mathrm{Ann\,}}

\newcommand{\tr}{\mathrm{tr}}

\newcommand{\Hom}{\mathrm{Hom}}

\newcommand{\Res}{\mathrm{Res}}
\newcommand{\End}{\mathrm{End}}

\newcommand{\cosoc}{\mathrm{cosoc}}

\newcommand{\SL}{\mathrm{SL}}

\newcommand{\m}{\mathfrak{m}}
\newcommand{\bT}{\mathbb{T}}

\newcommand{\Gal}{\mathrm{Gal}}

\newcommand\restr[2]{{% we make the whole thing an ordinary symbol
  \left.\kern-\nulldelimiterspace % automatically resize the bar with \right
  #1 % the function
  \right|_{#2} % this is the delimiter
  }}

\title[Modular character formulas and minimal lifts]{
Modular analogs of character formulas and minimal lifts of modular forms}

\author{Patrick B. Allen}
\address[Patrick Allen]{Department of Mathematics and Statistics, McGill University \\
Burnside Hall \\
805 Sherbrooke Street West \\
Montreal, Quebec H3A 0B9}
\email{patrick.allen@mcgill.ca}

\author{Preston Wake}
\address[Preston Wake]{Department of Mathematics, Michigan State University \\
619 Red Cedar Road \\
C212 Wells Hall \\
East Lansing, MI 48824}
\email{wakepres@msu.edu}

\date{\today}

\begin{document}

\begin{abstract}
If $f$ is a mod-$3$ eigenform of weight 2 and level $\Gamma_0(\ell^2)$ for a prime $\ell$ such that $\ell \equiv -1 \pmod{3}$, and $\ell$ is a vexing prime for $f$, we show that there is no obstruction to finding a minimal lift of $f$, but that there is an obstruction to finding a nonminimal lift.
The key new ingredient that we prove is a modular analog of the standard character formula for a cuspidal representation of $\GL_2(\F_\ell)$, an enhancement that allows us to easily compute the group cohomology of a $3$-adic lattice in such a representation. 
In fact, we provide a general framework for proving such modular analogs for a broader class of representations using results of Brou\'e and Puig in modular representation theory. We show that this class includes certain Deligne--Lusztig representations and 
representations coming from higher-depth supercuspidal representations of $\mathrm{GL}_2$.
\end{abstract}


\maketitle

\setcounter{tocdepth}{1}

\tableofcontents

%%%%%%%%%%%%%%%%%%%%%%%%%%%% Introduction %%%%%%%%%%%%%%%%%%%%%%%%%%%%

\section{Introduction}

In \cite{serresconj}, Serre formulated his celebrated conjecture that describes which mod-$p$ Galois representations arise from modular forms, and, moreover, the precise minimal weight, level, and character of such a form. Later, Serre \cite{serrenotes} found a counterexample to the character part of the conjecture for $p \in \{2,3\}$, arising from the stacky nature of modular curves for those primes. With this caveat, Serre's conjecture has been proven \cite{KW}. However, there are finer notions of minimality that come from Galois deformation theory and there is precisely one case where this Galois-theoretic minimality is not determined by the weight, level, and character, and this is the case of a vexing prime (so named in \cite{diamond} because of this property). One can ask whether the obstructions identified by Serre for $p \in \{2,3\}$ also obstruct minimal lifting in this sense. We show that this situation is the reverse of the situation considered by Serre: there is no obstruction to minimal lifting, but there is an obstruction to nonminimal lifting. The key to adapting Serre's construction to this setting is a formula for computing the group cohomology of cuspidal representations.

\subsection{Modular analogs of character formulas}
Many interesting irreducible characteristic-zero representations of finite groups are best described by character formulas. Typically, these formulas involve writing the character of the representation as a linear combination of characters of representations that are better understood (for instance, those that are induced from subgroups that are simpler). Since the category of characteristic-zero representations is semisimple, the character formula is essentially enough to completely understand the representation. However, if one wants to understand the modular theory (that is, over characteristic-$p$ fields), then the character formula can be insufficient. Indeed, for the application to the vexing-analog of Serre's conjecture, we need to compute the cohomology of the $p$-adic lattice in a cuspidal representation of $\GL_2(\F_\ell)$, where the character formula is well known, but not particularly useful when computing cohomology. We consider the question of finding \emph{modular analogs} of character formulas, and describe a situation where this question comes down to a question about multiplicity-freeness, and prove some results relating multiplicity-freeness in characteristic zero and characteristic $p$.

Before stating our general result, we introduce some notation that will be used throughout the paper: $p$ is a prime number, $G$ is a finite group, and $K$ is a $p$-adic local field that contains the $|G|$th-roots of unity, with valuation ring $\cO$, uniformizer $\varpi$, and residue field $k=\cO/\varpi\cO$. We also use the following naming convention: if $\tilde{M}$ is an $\cO[G]$-module, then $M=\tilde{M}\otimes_\cO K$ and $\ov{M}=\tilde{M} \otimes_\cO k$; note that a $K[G]$-module $M$ does not determine a $\cO[G]$-lattice $\tilde{M}$ in general, but if $M$ is a $K[G]$-module that has a unique (up to scaling) $\cO[G]$-lattice (for instance, if $M$ is one-dimensional), then we write
$\tilde{M}$ for that lattice.
We use the language of modular representation theory, which we review in Section \ref{subsec:review} below.

\begin{theorem}
\label{thm:main intro}
    Let $\pi$ be a simple $K[G]$-module that is in a nilpotent block with cyclic defect group. Let $e \in \cO[G]$ be the corresponding block idempotent. Suppose that there is a projective $\cO[G]$-module $\tilde{M}$, a subgroup $H \le G$, and a character $\tilde{\theta}: H \to \cO^\times$ such that there is an equality of virtual representations
    \begin{equation}
    \label{eq:character formula}
            \pi = M - \Ind_H^G\theta.
    \end{equation}

    For $s \in G$, let $H_s=H \cap s^{-1}Hs$ and let $\theta^s:H_s \to K^\times$ be $\theta^s(x)=\theta(sxs^{-1})$. Assume that $\Ind_H^G\theta$ is multiplicity free, and that there is no $s \in G$ such that the character $\theta^{-1}\theta^s:H_s \to K^\times$ is nontrivial and has $p$-power order. Then there is an $\cO$-lattice $\tilde{\pi}$ in $\pi$ such that there is an exact sequence of $\cO[G]$-modules
    \begin{equation}
    \label{eq:modular analog}
          0 \to \tilde{\pi} \to e\tilde{M} \to e \Ind_H^G\tilde{\theta} \to 0.
    \end{equation}
    In particular, for every subgroup $G' \le G$ and every $n \in \Z$ there is an isomorphism of Tate cohomology groups
    \begin{equation}
    \label{eq:tate cohomology}
            \hat H^n(G',\tilde{\pi}) \cong \hat H^{n-1}(G',{e}\Ind_H^G\tilde\theta).
    \end{equation}
\end{theorem}

There are many irreducible representations $\pi$ that have character formulas like \eqref{eq:character formula} (for instance, where $\tilde{M}$ is an induction from a prime-to-$p$-order subgroup). We think of \eqref{eq:modular analog} as a \emph{modular analog} of the character formula, in that it upgrades \eqref{eq:character formula} to a short exact sequence of lattices. The utility of this is illustrated by \eqref{eq:tate cohomology}: $H^{n-1}(G',\Ind_H^G\tilde\theta)$ can often be computed easily by applying Mackey's formula and Shapiro's lemma.  

We show that assumptions on $\pi$ are satisfied in many cases that are of interest in number theory. For instance, we show that if $G$ is a Chevalley group in characteristic $\ell$, $T$ is a maximal torus in $G$ such that the $p$-part of $T$ is a Sylow-$p$ subgroup of $G$, and $\theta$ is a character of $T$ that is general position modulo $p$, then the Deligne--Lusztig irreducible representation $\pi_\theta$ associated to $\theta$ is in a nilpotent block. This has interesting applications even for $G=\GL_2(\F_\ell)$.

\begin{ex}
\label{ex:GL2}
    Let $\ell$ be a prime number with $\ell \equiv -1 \pmod{p}$ and let $G=\GL_2(\F_\ell)$ and $T \subset G$ be a nonsplit torus isomorphic to $\F_{\ell^2}^\times$. Let $\theta: T \to K^\times$ be a character whose order is not a power of $p$ and let $\pi_\theta$ be the associated cuspidal representation of $G$.

    There is a well-known character formula for $\pi_\theta$:
    \[
    \pi_\theta = \Ind_{ZU}^G{\psi} - \Ind_T^G{\theta}.
    \]
    where $Z$ is the center of $G$, $U$ the upper-triangular unipotent subgroup, and $\psi$ is character that is nontrival on $U$ and agrees with $\theta$ on $Z$. In this case, Theorem \ref{thm:main intro} yields an exact sequence
\begin{equation}
\label{eq:GL2 sequence}
    0 \to \tilde{\pi}_\theta \to {e}\Ind_{ZU}^G\tilde\psi \to {e}\Ind_T^G\tilde\theta \to 0
\end{equation}
that can be used to compute the cohomology of $\tilde{\pi}_\theta$.
\end{ex}
\noindent We have similar applications to higher-depth supercuspidal types for $\GL_2(\Q_\ell)$.

The proof of Theorem \ref{thm:main intro} proceeds as follows. First, we reduce the problem to proving a bound on the multiplicities in the cosocle of $e\Ind_H^G \ov \theta$. The assumption that $\Ind_H^G\theta$ is multiplicity free implies that its endomorphism ring is commutative, and we use the additional assumptions on $\theta$ to show that the endomorphism ring of the reduction $\Ind_H^G \ov \theta$ is commutative. Finally, we use the Brou\'e--Puig structure theory of nilpotent blocks to deduce from this commutativity that the cosocle of $e\Ind_H^G \ov \theta$ is multiplicity free.

We note that Zhang \cite{zhang2024} had claimed a more general multiplicity-freeness result for reductions of induced representations, but his results are not true in the generality stated. We collaborated on a partial correction \cite{allen-wake-zhang-2025} to \cite{zhang2024}, but the correction is not general enough to cover the situation considered here.

\subsection{Modular forms and minimal lifts}
\label{subsec:minimal lifts}
The instance of Theorem \ref{thm:main intro} in Example \ref{ex:GL2} has an application to modular forms and Galois representations.
Let $\Gal_\Q=\Gal(\ov\Q/\Q)$ and let $\ov\rho:\Gal_\Q \to \GL_2(k)$ be an absolutely irreducible two-dimensional mod-$p$ Galois representation. Serre \cite{serresconj} formulated a precise conjecture about which $\ov \rho$ come from modular forms, and, when they do, what the minimal weight, level, and character of the form are. The philosophy for the conjecture is that the weight, level, and character should be determined locally and that there should be no global obstruction.
Later, he found that there was an error in the character part of his conjecture for $p \in \{2,3\}$ and gave a counterexample \cite[Contre-exemple, pg.~54]{serrenotes} for $p=3$ using the modular form \cite[Newform orbit 13.2.e.a]{lmfdb}. 
As Serre explains, the fact that $\SL_2(\Z)$ contains elements of order $2$ and $3$ introduces a global obstruction.

Theorem \ref{thm:main intro} can be used to study related questions. The minimal weight, level, and character are closely related to the notion of a \emph{minimally ramified lift} of $\ov \rho$. In fact, there is only one type of prime $\ell$, called \emph{vexing primes}, for which the minimally ramified lift is not determined by the weight, level, and character; these are the primes $\ell \equiv -1 \pmod p$ such that $\ov\rho|_{D_\ell}$ is irreducible and $\ov\rho|_{I_\ell}$ is reducible where $I_\ell \le D_\ell \le \Gal_\Q$ are an inertia group and decomposition group at $\ell$. For every vexing prime $\ell$ and  eigenform $f$ associated to $\ov \rho$, $f$ is supercuspidal at $\ell$ with character $\theta$ (more precisely, the local component at $\ell$ of the automorphic representation of $f$ is a supercuspidal representation whose Bernstein component is associated to some smooth character $\theta$ of $\Z_{\ell^2}^\times$). 
The minimally ramified (at $\ell$) lifts are the ones where the character $\theta$ has minimal order. The philosophy of Serre's conjecture would suggest that, unless $p \in \{2,3\}$, there should be no obstruction to finding lifts of either type---minimally ramified at $\ell$ or non-minimally ramified at $\ell$---and, indeed, one can show this. When $p \in \{2,3\}$, it is not clear whether the global obstruction as in Serre's counterexample plays a role. The following theorem shows that there is no obstruction to finding a minimal lift.

\begin{theorem}
    \label{thm:minimal lifts}
    Let $p=3$ and let $\ell$ be a prime with $\ell \equiv -1\pmod{3}$. Let $f$ be a newform of weight $2$ and level $\Gamma_0(\ell^2)$ that is depth-zero supercuspidal at $\ell$ with character $\theta: \F_{\ell^2}^\times \to \ov\Q_3^\times$ whose order is divisible by $3$. Let $\theta':\F_{\ell^2}^\times \to \ov\Q_3^\times$ be the character that has prime-to-$3$ order and that is equal to $\theta$ on the prime-to-$3$-order elements.
    
    Then there is a newform $g$ of weight $2$ and level $\Gamma_0(\ell^2)$ such that
    \begin{enumerate}
        \item The Hecke eigenvalues of $f$ and $g$ are congruent modulo a prime above $3$, and
        \item $g$ is supercuspidal at $\ell$ with character $\theta'$.
    \end{enumerate}
\end{theorem}
With the notation as in the theorem, let $\ov\rho$ be the mod-$3$ Galois representation associated to $f$ and let $\rho_g$ be the $3$-adic Galois representation associated to $g$. Then (1) says that $\rho_g$ is a lift of $\ov\rho$ and (2) says that it is a minimal lift. Hence, there no obstruction to finding a minimal lift. However, the following example shows that the roles of $f$ and $g$ in the theorem cannot be reveresed: there \emph{is} an obstruction to finding non-minimal lifts.

\begin{ex}
\label{ex:3 11}
    Let $p=3$ and $\ell=11$. Consider the newform $g$ of weight $2$ and level $\Gamma_0(121)$ with $q$-expansion
    \[
    g = q - q^3 - 2   q^4 - 3   q^5 - 2   q^9 + 2   q^{12} + \dots
    \]
    \cite[Newform orbit 121.2.a.b]{lmfdb}. Using the Loeffler--Weinstein algorithm \cite{LW2012} in Sage \cite{sagemath} reveals that $g$ is depth-zero supercuspidal at $11$ with character $\theta$ of $\F_{121}^\times$ that has order $4$. In particular, the reduction $\ov\theta$ modulo $3$ also has order $4$ and hence cannot factor through $\F_{11}^\times$, so it is in general position.

    Let $\rho_g$ and $\ov \rho$ be the $3$-adic and mod-$3$ representations of $g$, respectively. The above calculation shows that $11$ is a vexing prime for $\ov \rho$ and that $\rho_g$ is a minimally ramified lift of $\ov \rho$. Examining the list of newforms of weight $2$ and level $\Gamma_1(121)$ in \cite{lmfdb} quickly reveals that there are no other forms there congruent to $g$. In particular, in this case, there are no lifts of $\ov\rho$ of weight $2$ and level $\Gamma_1(121)$ that are \emph{not} minimally ramified at $11$.
\end{ex}

In Serre's counterexample, there is an obstruction to finding a minimal lift, but no obstruction to finding a non-minimal lift. In the vexing situation we consider here, it is the reverse. This reversal is accounted for by the dimension shift in \eqref{eq:tate cohomology} of Theorem \ref{thm:main intro}.

Theorem \ref{thm:minimal lifts} raises many more questions. What happens for higher-depth supercuspidals? What happens if there are both vexing primes and principal series primes? In an obstructed situation, can the obstruction be made explicit? This circle of ideas will be explored more fully in our forthcoming work. More ambitiously, one can consider analogous questions for automorphic forms on reductive groups $G$ other than $G=\GL_2$. In this case, the analog of the cuspidal representations $\pi_\theta$ are the Deligne--Lusztig representations considered in Section \ref{sec:DL}.

\subsection{Outline of the paper}
We give the proof of Theorem \ref{thm:main intro} in Section \ref{sec:proof}. In Section \ref{sec:DL}, we use Deligne--Lusztig theory to show that our results apply to many representations that are of interest in number theory. In Section \ref{sec:GL2}, we apply our results to the case of $\GL_2$, and in Section \ref{sec:min lifts}, we prove our application to modular forms.

\subsection{Acknowledgements}
P.~W.~thanks Bob Kottwitz, Jackie Lang and Rob Pollack for helpful conversations
about modular representation theory and Deligne--Lusztig theory. Both authors thank Robin Zhang for clarifying discussions about \cite{zhang2024} and for collaborating on \cite{allen-wake-zhang-2025}.

P.~A.~acknowledges the support of the Natural Sciences and Engineering Research Council of Canada (NSERC), [funding reference number RGPIN-2020-05915].
P.~A.~a été financée par le Conseil de recherches en sciences naturelles et en génie du Canada (CRSNG), [numéro de référence RGPIN-2020-05915].
P.~W.~is supported by National Science Foundation CAREER Grant No.~DMS-2337830.

\section{Modular analogs}
\label{sec:proof}
In this section, we prove Theorem \ref{thm:main intro}. We first prove a general result relating modular analogs of character formulas to multiplicities of cosocles. Then we use a result of Brou\'e and Puig to relate multiplicity freeness of cosocles to commutativity of endomorphism rings. We then complete the proof using a commutativity result from \cite{allen-wake-zhang-2025}.

\subsection{Modular analogs and multiplicities} The following elementary lemma reduces the problem to a question about multiplicities. Recall the notational convention that if $\tilde{M}$ is an $\cO$-module, then $M$ denotes $\tilde{M} \otimes_\cO K$ and $\ov{M}$ and denotes $\tilde{M} \otimes_\cO k$.

\begin{lemma}
\label{lem:general modular analog} Let $\tilde{R}$ be a $\cO$-algebra that is free of finite rank as an $\cO$-module, and assume that $R=\tilde{R}\otimes_\cO K$ is semisimple. 
Let $\tilde{M}$ and $\tilde{M}'$ be finitely generated $\tilde{R}$-modules that are $\cO$-free. Assume that there is a surjective $R$-module homomorphism $M \onto M'$ and let $\pi$ denote its kernel. 

Assume that $\tilde{M}$ is projective and that, for every simple $\ov{R}$-module $\eta$,
\begin{equation}
\label{eq:mults}
        \dim_{k} \Hom_{\ov R}(\ov M,\eta) \ge \dim_k \Hom_{\ov R}(\ov {M'},\eta).
\end{equation}
    Then there is an $\cO$-lattice $\tilde{\pi}$ in $\pi$ such that there is a short exact sequence of $\tilde{R}$-modules
    \[
    0 \to \tilde{\pi}\to \tilde{M} \to \tilde{M}' \to 0.
    \]
\end{lemma}
\begin{proof}
    The assumption \eqref{eq:mults} implies that there is a surjective $\ov R$-module homomorphism $\ov{f}:\cosoc(\ov M) \onto \cosoc(\ov{M'})$. Since $\tilde{M}$ is projective, there is $\tilde{R}$-module homomorphism $f:\tilde{M} \to \tilde{M}'$ such that $\cosoc(f\otimes_\cO k)=\ov f$. By Nakayama's lemma, $f$ is surjective. Since $\pi=M-M'$, the kernel of $f$ is a lattice in $\pi$.
\end{proof}

\subsection{Review of modular representation theory}
\label{subsec:review}
We will make use of Brauer's theory of blocks. This theory is classical, but may not be familiar to the audience of this paper (especially those interested in the number theory applications). We encourage the interested reader to consult \cite[Chapter 15]{isaacs}, \cite{alperin}, or \cite{curtis-reiner} (see also \cite{bonnafe} for a focus on $\SL_2$). We review some of the relevant ideas:
\begin{itemize}
    \item Simple $k[G]$-modules correspond to irreducible \emph{Brauer characters} in a similar way that simple $K[G]$-modules correspond to ordinary characters. Brauer characters essentially restrictions of ordinary characters to elements of prime-to-$p$ order. (See \cite[Theorem 15.6, pg.~265]{isaacs}.)
    \item The irreducible representations $\rho$ of $G$ over $K$ are grouped together into \emph{blocks} according to the simple $k[G]$-modules that appear as Jordan--Holder factors of the reduction of an $\cO$-lattice in $\rho$. (See \cite[Definition 15.17 pg.~271]{isaacs}.)
    \item Blocks are characterized by integrality of idempotents. Associated to $\rho$ is the central idempotent
    \[
    e_\rho = \frac{\dim_K(\rho)}{\#G} \sum_{g \in G} \tr(\rho(g))[g^{-1}]
    \]
    in $K[G]$; these are mutually orthogonal. Due to the $\#G$ in the denominator, this will not typically be an element of $\cO[G]$. A set $B$ of irreducible representations is a union of blocks if and only of $e_B=\sum_{\rho \in B} e_\rho$ is in $\cO[G]$. (See \cite[Theorems 15.26, 15.27, 15.30, pg.~275]{isaacs}.)
    \item Associated to a block $B$ is a conjugacy class of $p$-groups $D$ called the defect groups of the block. The $p$-part of the index of $D$ in $G$ is given by
    \begin{equation}
    \label{eq:defect order}
    [G:D]_p = \max_a \{p^a  \ : \ p^a \text{ divides }\dim_K(\rho) \text{ for all }\rho\in B\}.
    \end{equation}
    (See \cite[Theorem 15.41, pg.~280]{isaacs}.)
    \item Brauer related blocks of $G$ to blocks of certain subgroups using the Brauer homomorphism. If $P \le G$ is a $p$-subgroup, let $C_G(P)$ be its centralizer. The $k$-linear restriction map
    \[
    \mathrm{Br}_P: Z(k[G]) \to Z(k[C_G(P)])
    \]
    is a ring homomorphism, called the \emph{Brauer homomorphism} (see \cite[Lemma 15.32, pg.~278]{isaacs}). If $B$ is a block of $G$ and $B'$ is a block of $C_G(P)$, then $(P,B')$ is called a \emph{$B$-Brauer pair} if $e_{B'} \mathrm{Br}_P(e_B)=e_{B'}$. Then $G$ acts on the set of $B$-Brauer pairs by conjugation, and the stabilizer of $(P,B')$ is denoted $N_G(P,B')$; it plays a role in block theory similar to that of $N_G(P)$ in group theory. Brauer's Second Main theorem (\cite[Theorems 15.48, pg.~284]{isaacs}) is an important tool for identifying $B$-Brauer pairs.
\end{itemize}

\subsection{Nilpotent blocks} 
Recall that $G$ is called \emph{$p$-nilpotent} if $G=P \ltimes H$ for a Sylow-$p$ subgroup $P$ and a normal subgroup $H$. It is a classical theorem of Frobenius that $G$ is $p$-nilpotent if and only if $N_G(P)/C_G(P)$ is a $p$-group for every $p$-subgroup $P \le G$. Based on this, Brou\'e and Puig have made the following definition \cite{BP1980}: a block $B$ for $G$ is \emph{nilpotent} if $N_G(P,f)/C_G(P)$ is a $p$-group for every $B$-Brauer pair $(P,f)$.

They show that nilpotent blocks are particularly simple: for instance, they only include a single simple $k[G]$-module (see \cite[Theorem 1.2]{BP1980}). Moreover, they show the following in \cite[pg.~120]{BP1980}.

\begin{theorem}[Brou\'e--Puig]
\label{thm:BP}
    Let $B$ be a nilpotent block that has an abelian defect group $D$. Then ${e}_B\cO[G]$ is Morita equivalent to $\cO[D]$.
\end{theorem}

The statement of this theorem may seem abstract, but, in practice, the map $\tilde{e}_B\cO[G] \to\mathrm{M}_n(\cO[D])$ can be made completely explicit and it is not difficult to see why it is an isomorphism; we illustrate this in a particular case in Remark \ref{remark:explicit nilpotent block} below.

Brou\'e and Puig \cite[pg.~118]{BP1980} also give the following criteria for $B$ to be nilpotent, the first of which is clear from the definition and Frobenius's theorem, and the second of which is based on \cite[Proposition 4.21]{AB1979}.

\begin{theorem}[Alperin--Brou\'e]
\label{thm:nilpotent blocks examples} Let $B \in \mathrm{Bl}(G)$.
\hfill
\begin{enumerate}
    \item If $G$ is $p$-nilpotent, then $B$ is nilpotent.
    \item If there is a $B$-Brauer pair $(D,f)$ where $D$ is a defect group of $B$ such that $N_G(D,f)=C_G(D)$ and $D$ is abelian, then $B$ is nilpotent.
\end{enumerate}
\end{theorem}

\subsection{Multiplicity-freeness over monogenic local $k$-algebras}
According to Theorem \ref{thm:BP}, if $B$ is nilpotent block that has a cyclic defect group $D$, then ${e}_Bk[G]$ is Morita equivalent to $k[D] \cong k[x]/(x^{\#D})$. Modules over this kind of ring are, of course, particularly simple.

\begin{lemma}
\label{lem: monogenic}
Let $S=k[x]/(x^m)$ for some $m>0$, and let $R$ be a $k$-algebra such that there is a Morita equivalence of categories
\[
F: S\text{-mod} \to R\text{-mod}. 
\]
and let $\pi=F(k)$. Let $M$ be a finitely-generated nonzero $R$-module. Then $\cosoc(M)\cong\pi$ if and only if $\End_R(M)$ is commutative.
\end{lemma}
\begin{proof}
By the equivalence, it suffices to consider the case $R=S$. In that case, if $\cosoc(M)\cong k$, then $M$ is cyclic by Nakayama's lemma and, since $S$ is commutative, $\End_S(M) \cong S/\mathrm{Ann}_S(M)$ is commutative. Conversely, if $\End_S(M)$ is commutative, then it is easy to see that $M$ must be indecomposable. Indeed, if $M=A \oplus B$ for nonzero $A$ and $B$, then, using the fact that $\Hom_S(A,B)$ is nonzero, it is easy to produce non-commutating elements of $\End_S(M)$. Since every finitely-generated $S$-module is a direct sum of cyclic modules, this implies that $M$ is cyclic, and hence that $\cosoc(M)\cong k$.
\end{proof}

\subsection{Commutativity of endomorphisms of induced representations}
\label{subsec:commutative}
With the notation and assumptions as in Theorem \ref{thm:main intro}, \cite[Proposition 2.2]{allen-wake-zhang-2025} implies that $\End_{k[G]}(\Ind_H^G \ov \theta)$ is commutative. We recall the idea of the proof here. The assumption that $\Ind_H^G\theta$ is multiplicity free implies that $\End_{K[G]}(\Ind_H^G\theta)$ is commutative, and hence that the subring $\End_{\cO[G]}(\Ind_H^G \tilde{\theta})$ is commutative. A Mackey-theory argument shows that, under the additional assumptions on $\theta$, the natural map
\[
\End_{\cO[G]}(\Ind_H^G \tilde{\theta}) \to \End_{k[G]}(\Ind_H^G\ov \theta)
\]
is surjective. This implies that $\End_{k[G]}(\Ind_H^G \ov \theta)$ is commutative.

\subsection{Proof of Theorem \ref{thm:main intro}}
As was just recalled in Section \ref{subsec:commutative}, the assumptions imply that $\End_{k[G]}(\Ind_H^G \ov \theta)$ is commutative. It follows that $\End_{ek[G]}(e\Ind_H^G \ov \theta)$ is commutative. By Theorem \ref{thm:BP} of Brou\'e--Puig, $ek[G]$ is Morita-equivalent to a ring of the form $k[x]/(x^m)$, so we may apply Lemma \ref{lem: monogenic} to conclude $\cosoc(e\Ind_H^G \ov \theta)\cong \ov\pi$. Hence, we may apply Lemma \ref{lem:general modular analog} with $\tilde{R}=e\cO[G]$ and $\tilde{M}' = e\Ind_H^G \tilde{\theta}$ to obtain the desired exact sequence \eqref{eq:modular analog}. The formula \eqref{eq:tate cohomology} follows from the fact that $e\tilde{M}$, being projective, is cohomologically trivial.\hfill\qed

\section{Nilpotent blocks in Deligne--Lusztig theory}
\label{sec:DL}
We use Deligne--Lusztig theory to show that certain Deligne--Lusztig characters are in nilpotent blocks with cyclic defect group. This allows us to apply the results of the previous section to these characters.

\subsection{Deligne--Lusztig theory} We summarize some of the main results of Deligne--Lusztig theory (see \cite{DL1976} and \cite{serre1977}). Let $\mathcal{G}$ be a reductive group over finite field $F$ of characteristic $\ell$ with a maximal torus $\mathcal{T} \subset \mathcal{G}$. Let $G=\mathcal{G}(F)$ and $T=\mathcal{G}(F)$ and let $W=N_G(T)/T$ be the Weyl group of $T$.

Let $\theta:T \to K^\times$ be a character. For $w \in W$, let $\theta^w(t)=\theta(w t w^{-1})$. We say that $\theta$ is \emph{in general position} if $\theta^w \ne \theta$ for all nontrivial $w \in W$.
The following theorem summarizes the parts of Deligne--Lusztig theory that we will need.

\begin{theorem}[Deligne--Lusztig]
\label{thm:DL summary}
For each character $\theta: T \to K^\times$, there is a class function
\[
R_\theta: G \to K
\]
that is an integer linear combination of trace-characters of irreducible representations of $G$ and satisfies:
\begin{itemize}
    \item If $u \in G$ is unipotent, then $R_\theta(u)$ is an integer, independent of $\theta$, and is denoted $Q_{T,G}(u)$. The function $Q_{T,G}$ and is called the \emph{Green function} of $G$ and $T$.
    \item If $x \in G$ is written as $x=su=us$ for $u \in G$ unipotent and $s \in G$ semisimple, then
\begin{equation}
\label{eq:D-L formula}
    R_\theta(su) =  \sum_{\substack{g\in G/C_G(s) \\ g^{-1}sg \in T}} Q_{gTg^{-1},C_G(s)}(u) \theta(g^{-1}sg).
\end{equation}
    \item If $\theta$ is in general position, then there is a $\epsilon_{G,T} \in \{\pm 1\}$, independent of $\theta$, such that $\epsilon_{G,T}R_\theta$ is the trace character of an irreducible representation $\pi_\theta$ of $G$ of dimension $\frac{\#G}{\#T\#L}$ where $L\le G$ is an $\ell$-Sylow subgroup. In particular, if $G$ itself is a torus and $T=G$, then $\pi_\theta=\theta$ and $Q_{T,G}(1)=1$.
    \item If $\theta$ is in general position and $\mathcal{T}$ is not contained in a proper parabolic subgroup of $\mathcal{G}$, then $\pi_\theta$ is cuspidal. 
\end{itemize}
\end{theorem}
\begin{proof}
The first two points follow from \cite[Theorem 4.2]{DL1976} (the integrality of $R_\theta(u)$ follows from \cite[Proposition 3.3]{DL1976}). In the third point, the irreducibility follows from \cite[Theorem 6.8]{DL1976}, and the dimension formula from \cite[Theorem 7.1]{DL1976}. The final point is \cite[Theorem 8.3]{DL1976}.
\end{proof}

\subsection{Nilpotent blocks of Deligne--Lusztig characters} Let $\ell$, $G$, and $T$ be as in the previous section. In this section, we'll consider the case where the $p$-Sylow subgroup of $T$ is a $p$-Sylow subgroup of $G$ that contains a regular semisimple element.

\begin{ex}
    Let $F=\F_q$ for $q$ a power of $\ell$ and let $n$ be the multiplicative order of $q$ modulo $p$. Let $\mathcal{G}=\GL_n$ and let $\mathcal{T} \le \mathcal{G}$ be a torus isomorphic to $\Res_{\F_{q^n}/\F_q}(\G_m)$. Then $T=\F_{q^n}^\times$ is a cyclic group whose generators are regular semisimple elements of $G$ and whose $p$-Sylow subgroup is a $p$-Sylow subgroup of $G$.
\end{ex}

For a character $\theta: T \to \cO^\times$, let $\ov\theta:T \to k^\times$ be the reduction modulo $\varpi\cO$.
\begin{theorem}
\label{thm:D-L and nilpotent blocks}
Let $\theta:T \to \cO^\times$ be a character. Assume that
\begin{itemize}
    \item The $p$-Sylow subgroup $P$ of $T$ is a $p$-Sylow subgroup of $G$,
    \item $P$ contains a regular semisimple element of $G$, and
    \item $\ov \theta$ is in general position.
\end{itemize}
        Let $\Xi$ be the set of $p$-power-order $\cO$-valued characters of $T$. Then $\theta\xi$ is in general position for all $\xi \in \Xi$ and the collection $\{\pi_{\theta \xi} : \xi \in \Xi\}$ of irreducible representations forms a nilpotent block with defect group $P$.
\end{theorem}

\begin{proof}
Since the characters $\theta\xi$ all reduce to $\ov{\theta}$, they are all in general position. By Theorem \ref{thm:DL summary}, there is an irreducible representation $\pi_{\theta\xi}$ of $G$ whose trace-character $\chi_{\theta\xi}$ is $\epsilon_{G,T} R_{\theta\xi}$ for a sign $\epsilon_{G,T} \in \{\pm1\}$ that is independent of $\xi$ and $\theta$, and whose dimension $\chi_{\theta\xi}(1)$ is $\frac{\#G}{\#T \#L}$, which is prime-to-$p$ (by the assumption that $P\le T$ is a $p$-Sylow of $G$) and independent of $\xi$ and $\theta$.

Let 
\[
e_{\theta\xi} = \frac{\chi_{\theta\xi}(1)}{\#G} \sum_{x \in G} \chi_{\theta\xi}(x) [x^{-1}]
\]
be the idempotent in $K[G]$ associated to $\pi_{\theta\xi}$, and let $\tilde{e}= \sum_{\xi \in \Xi} e_{\theta\xi}$. We claim that $\tilde{e} \in \cO[G]$. To see this, note that, by the above-summarized properties of $\chi_{\theta\xi}$, the idempotent $\tilde{e}$ can be written as 
\[
\tilde{e} = \epsilon_{G,T}\frac{\chi_\theta(1)}{\# G} \sum_{x \in G} \tilde{R}_\theta(x) [x^{-1}]
\]
where $\tilde{R}_\theta = \sum_{\xi \in \Xi} R_{\theta\xi}$. It's clear then that $\tilde{e} \in \frac{1}{\#G}\cO[G]=\frac{1}{\#P} \cO[G]$, and to show that $\tilde{e}$ is in $\cO[G]$, it is enough to show that $\tilde{R}_\theta(x) \in \#P \cdot\cO$ for all $x \in G$. Writing $x=su=us$ with $s$ semisimple and $u$ unipotent, the formula \eqref{eq:D-L formula} for $R_{\theta\xi}$ implies
\begin{align*}
    \tilde{R}_\theta(x) &=  \sum_{g } Q_{gTg^{-1},C_G(s)}(u) \theta(g^{-1}sg)\left( \sum_{\xi \in \Xi}\xi(g^{-1}sg)\right) \\
    &= \begin{cases}
    \displaystyle \#\Xi \cdot \sum_{g } Q_{gTg^{-1},C_G(s)}(u) \theta(g^{-1}sg) & \text{if $s$ has prime-to-$p$ order} \\
    0 & \text{otherwise}
\end{cases},
\end{align*}
where the second equality uses the usual character-sum formula. 
Since $\#\Xi=\#P$, and the values of $\theta$ and the Green functions $Q$ are integral, this shows that $\tilde{R}_\theta(x)~\in~\#P \cdot \cO$ for all $x \in G$, as required.

As we recalled in Section \ref{subsec:review}, the fact that $\tilde{e}$ is in $\cO[G]$ implies that $\{\pi_{\theta\chi} : \xi \in \Xi\}$ is a union of blocks. But it is clear from the character formula \eqref{eq:D-L formula} that the associated Brauer characters $\hat{\chi}_{\theta\xi}$ are independent of $\xi$, so the representations $\pi_{\theta\xi}$ are contained in the same block by definition. This shows that $B=\{\pi_{\theta\xi} : \xi \in \Xi\}$ is a block with ${e}_B=\tilde{e}$. Since the representations $\pi_{\theta\xi}$ in the block have prime-to-$p$ dimension, \eqref{eq:defect order} implies that every defect group of $B$ is a $p$-Sylow subgroup of $G$. Since the defect groups of a block are a conjugacy class and every $p$-Sylow subgroup of $G$ is conjugate to $P$, this implies that $P$ is a defect group of $B$.

To show that $B$ is a nilpotent block, we will find a $B$-Brauer pair and apply the Alperin--Brou\'e criterion to it. The group $C_G(P)$ is $T$; let $f$ be the block of $k[T]$ that contains the character $\theta$. We claim that $(P,f)$ is an $B$-Brauer pair. To check this, we use Brauer's Second Main Theorem.
%with $x$ a regular element in $P$, $\rho=\pi_\theta$ and $\varphi=\hat{\ov\theta}$. 
Write $T=P \times S$ where $S$ is the prime-to-$p$ subgroup of $T$. The simple $k[T]$-modules are in bijection with $\hat{S}=\Hom(S,k^\times)=\Hom(T,k^\times)$, and, for $\eta \in \hat{S}$, its Brauer character $\hat\eta:S \to \cO$ is the composition $S \xrightarrow{\eta} k^\times \to \cO^\times$ of $\eta$ with the canonical splitting. Let $x$ be a regular element of $P$. By Brauer's Second Main Theorem, to show that $(P,f)$ is a $B$-Brauer pair, we need only show $d_{\ov\theta} \ne 0$, where, for $\eta \in \hat S$, $d_\eta$ is the generalized decomposition number defined by the equation
\begin{equation}
\label{eq:d_eta}
    \epsilon_{G,T}R_\theta(xy) = \sum_{\eta \in \hat{S}} d_\eta \hat\eta(y)
\end{equation}
for all $y \in S$.
Since $xy$ is regular, $C_G(xy)=T$ and $\{g \in G\colon g^{-1}xyg \in T\}=N_G(T)$, so \eqref{eq:D-L formula} reads
\begin{equation}
\label{eq:R(xy)}
    R_\theta(xy)=\sum_{w \in W(T)} \theta^w(xy).
\end{equation}
Combining \eqref{eq:d_eta} and $\eqref{eq:R(xy)}$ and reducing modulo $\varpi \cO$ yields
\[
\sum_{\eta \in \hat{S}} \ov d_\eta \eta(y) = \epsilon_{G,T}  \sum_{w \in W(T)} \ov{\theta^w}(y),
\]
where $\ov d_\eta$ is the image of $d_\eta$ in $k$.
Since the characters $\ov{\theta^w}$ are distinct by assumption, linear independence of characters implies that $\ov d_{\ov \theta}=\epsilon_{G,T} \ne 0$, so $(P,f)$ is an $e$-Brauer pair.

Finally, $N_G(P,f)/C_G(P)$ is the set of elements in $W(T)$ that preserve $\ov{\theta}$, which is trivial by the assumption that $\ov{\theta}$ is in general position. Hence the hypotheses of Theorem \ref{thm:nilpotent blocks examples}(2) are satisfied and $e$ is a nilpotent block.
\end{proof}

\begin{remark}
\label{remark:explicit nilpotent block}
For the sake of being self-contained, we make the isomorphism $\tilde{e}_B\cO[G] \cong \mathrm{M}_n(\cO[P])$ of Theorem \ref{thm:BP} completely explicit in this case, under the additional assumption that $P$ is cyclic, and describe the steps needed to show it is an isomorphism. First, there is the $\cO$-algebra homomorphism $\iota: \cO[P] \to Z(\tilde{e}_B\cO[G])$ given by $\iota(g)=\sum_{\xi \in \Xi} \xi(g)e_{\theta\xi}$; a priori $\iota(g) \in Z(\tilde{e}_BK[G])$, but a similar argument to the proof of the theorem shows that this element is integral. 
If $x \in P$ is a generator, then one can calculate that
\[
\iota(x-1)^{|P|-1} \not \equiv 0 \pmod{\varpi \cO[G]},
\]
and it follows that $\iota$ is split-injective as $\cO$-modules.
If $\mathcal{P}$ denotes a projective $\cO[G]$-module such that $\mathcal{P}\otimes_\cO k$ is a projective cover of $\ov\pi_\theta$, then it is easy to see that $\tilde{e}_B\cO[G] \cong \mathcal{P}^n$ for some $n$, so that the map $ Z(\tilde{e}_B\cO[G]) \to \End_{\cO[G]}(\mathcal{P})$ is split-injective. Finally, one can show that 
\[
\mathcal{P}\otimes_\cO K \cong \bigoplus_{\xi \in \Xi} \pi_{\xi\theta}
\]
so that the $\cO$-rank of $\End_{\cO[G]}(\mathcal{P})$ is $\#\Xi=\#P$. Hence the composition of split-injections
\[
\cO[P] \to Z(\tilde{e}_B\cO[G]) \to \End_{\cO[G]}(\mathcal{P})
\]
is an isomorphism. This shows that 
\[
\tilde{e}_B\cO[G] \cong \End_{\tilde{e}_B\cO[G]}(\mathcal{P}^n) \cong \mathrm{M}_n(\End_{\cO[G]}(\mathcal{P})) \cong \mathrm{M}_n(\cO[P]).
\]
\end{remark}

\section{Application to $\GL_2(\F_\ell)$ and to supercuspidal representations}
\label{sec:GL2}

We illustrate how the results of the previous section may be applied in the simple example where $\mathcal{G}=\GL_2$. We also apply our results to the case of a $p$-nilpotent group, and explain how this can be used to the study of higher-depth supercuspidal representations of $\ell$-adic groups.

\subsection{Cuspidal representations of $\GL_2(\F_\ell)$}
 Let $\ell$ be a prime with $\ell \equiv -1 \pmod{p}$, and let $G=\GL_2(\F_\ell)$ and $T\cong \F_{\ell^2}^\times$ a non-split maximal torus in $G$. Let $Z \subset G$ be the center (the scalar matrices) and let $U \subset G$ be the upper-triangular unipotent matrices.

Let $\ov\theta: \F_{\ell^2}^\times \to k^\times$ be a character that does not factor through the norm $\F_{\ell^2}^\times \to \F_\ell^\times$, and let $\tilde\theta: \F_{\ell^2}^\times \to \cO^\times$ be a lift of $\ov\theta$. Let $\tilde\pi_\theta$ be an $\cO$-lattice in the cuspidal representation of $G$ associated to $\theta$ (this lattice is unique up to homothety because the $\tilde\pi_\theta \otimes_\cO k$ is irreducible). Let $\tilde\psi:ZU \to \cO^\times$ be a character that agrees with $\tilde\theta$ on $Z$ and is non-trivial on $U$. Let $B$ be the block for $G$ that contains $\pi_\theta$ and let ${e}_B \in \cO[G]$ be the associated idempotent.

\begin{corollary}
With the above notation, there is an exact sequence of $\cO[G]$-modules
\[
0 \to \tilde\pi_\theta \to {e}_B\,\Ind_{ZU}^G\tilde\psi \to {e}_B\,\Ind_T^G\tilde\theta \to 0.
\]
\end{corollary}
\begin{proof}
We verify the hypothesis of Theorem \ref{thm:main intro} with $\pi=\pi_\theta$, $H=T$, and $\tilde{M}=\Ind_{ZU}^G\tilde\psi$:
\begin{itemize}
    \item Equation \eqref{eq:character formula} is the well-known character formula $\pi_\theta = \Ind_{ZU}^G \psi - \Ind_T^G \theta$ (see \cite[Theorem, Section 6.4, pg.~47]{BH2006}, for instance),
    \item  $\pi_\theta$ is in a nilpotent block with cyclic defect group by Theorem \ref{thm:D-L and nilpotent blocks},
    \item $ZU$ has prime-to-$p$ order, so Maschke's theorem implies that $\tilde{\psi}$ is a projective $\cO[ZU]$-module; since induction preserves projectivity, $\tilde{M}$ is a projective $\cO[G]$-module;
    \item For all $s \in G$, $\theta^{-1}\theta^s: T_s \to K^\times$ is either trivial or does not have $p$-power order. Indeed, if $s$ is in $T$, then $\theta^{-1}\theta^s$ is trivial, and if $s$ is not in the normalizer of $T$, then $T_s=Z$ is a prime-to-$p$ group, so the order of $\theta^{-1}\theta^s$ is prime to $p$. If $s$ is in the normalizer of $T$ but not in $T$, then the assumption that $\ov\theta$ is in general position implies that $\ov\theta \ne \ov\theta^s$, so $\theta^{-1}\theta^s$ does not have $p$-power order.\qedhere
\end{itemize}
\end{proof}

\subsection{Representations of $p$-nilpotent groups and applications to supercuspidal representations}
Cuspidal representations of $\GL_2(\F_\ell)$ give rise to depth-zero supercuspidal representations of $\GL_2(\Q_\ell)$. To study higher-depth supercuspidal representations, instead of using Deligne--Lusztig theory, we can rely on the fact that relevant groups are $p$-nilpotent. We have the following general result for $p$-nilpotent groups.
\begin{theorem}
\label{thm:p-nilpotent case}
    Let $G$ be a finite group with a prime-to-$p$-order normal subgroup $N$ and cyclic $p$-Sylow subgroup $P$ such that $G=N \rtimes P$. Let $H \le N$ be a subgroup and $\tilde{\rho}:G \to \GL_n(\cO)$ be a representation of $G$ such that $\rho|_H$ is irreducible and let $B$ be the block of $G$ that contains $\rho$. Let $L \le G$ be a subgroup and let $\tilde{\lambda}:L \to \cO^\times$ be a character of prime-to-$p$ order and assume that $\Ind_L^G\lambda$ is multiplicity-free.

    Then there is a surjective $\cO[G]$-module homomorphism
    \[
    e_B\Ind_H^G(\tilde{\rho}|_H) \to e_B\Ind_L^G\tilde{\lambda}.
    \]
\end{theorem}
\begin{proof}
    Since $G$ is $p$-nilpotent, $B$ is a nilpotent block, and the defect groups of $B$ are cyclic (since they are subgroups of $P$). Just as in Section \ref{subsec:commutative} above, \cite[Proposition 2.2]{allen-wake-zhang-2025} implies that $\End_{k[G]}(\Ind_L^G\ov\lambda)$ is commutative. By Theorem \ref{thm:BP} and Lemma \ref{lem: monogenic}, $e_B\Ind_L^G\ov\lambda$ is either zero (in which case there is nothing to prove) or has cosocle $\ov \rho$.
    Since $\#H$ is prime-to-$p$, ${e}_B\Ind_H^G(\tilde\rho|_H)$ is a projective $\cO[G]$-module, and it is nonzero since $\Ind_H^G(\rho|_H)$ contains $\rho$. The result now follows from Lemma \ref{lem:general modular analog}.
\end{proof}

This theorem is potentially useful when understanding the modular representation theory of positive-depth supercuspidal representations of $\GL_n$, as constructed by Howe \cite{howe1977}. In the following example, we consider the case of $\GL_2$, following Bushnell--Henniart \cite{BH2006}.

\begin{ex}
\label{ex:higher depth}
    Let $\ell$ be a prime with $\ell \equiv -1 \pmod{p}$ and let:
    \begin{itemize}
        \item $H$ be an $\ell$-group with cyclic center $Z$, such that $V=H/Z$ is a two-dimensional $\F_\ell$-vector space,
        \item $\theta: Z \into K^\times$ be a faithful character,
        \item $A$ be a cyclic group of order $\ell+1$ that acts on $H$ in such a way that it acts on $Z$ trivially, and, for all nontrivial $a \in A$, the $a$-fixed part $V^a$ of $V$ is zero.
    \end{itemize}
Let $G=H \rtimes A$. In this situation, there is a unique irreducible representation $\eta$ of $H$ such that $\eta|_Z$ contains $\theta$ \cite[Lemma 2 of Section 16.4, pg.~114]{BH2006}. This representation therefore extends to a representation of $G$, but this extension is not unique. In \cite[Lemma of Section 22.2, pg.~135]{BH2006}, it is shown that there is a unique extension $\eta_1$ of $\eta$ to $G$ such that $\tr(\eta_1(a))=-1$ for all nontrivial $a \in A$, and the construction of this $\eta_1$ is a key step in constructing supercuspidal representations of $\GL_2(\Q_p)$ of positive even depth. In the course of the proof of \emph{loc.~cit.}, the representation $\eta_1$ is described by a character formula
\[
\eta_1 = \Ind_H^G \eta-\Ind_{ZA}^G\theta
\]
and it is shown that $\Ind_{ZA}^G\theta$ is multiplicity-free. Using Theorem \ref{thm:p-nilpotent case}, we can make an analogous statement to this character formula in the context of $\cO[G]$-modules. Indeed, writing $A=P \times A'$, where $P$ is the $p$-Sylow subgroup of $A$, we can apply Theorem \ref{thm:p-nilpotent case} with $N=H \rtimes A'$, $\rho={\eta}_1$, $L=ZA$, and $\lambda=\theta$ to obtain an exact sequence of $\cO[G]$-modules
\begin{equation}
\label{eq:eta1}
   0 \to \tilde\eta_1 \to {e}_B \Ind_H^G \tilde\eta \to {e}_B\Ind_{ZA}^G\tilde\theta \to 0.
\end{equation}
\end{ex}

\begin{ex}
We now explain how Example \ref{ex:higher depth} relates to even-depth supercuspidal representations of $\GL_2(F)$, where $F$ is an $\ell$-adic local field. To simplify the discussion, we only consider the case where $\ell>2$ and $F=\Q_\ell$. Let $\Q_{\ell^2}=W(\F_{\ell^2})[1/\ell]$ be the unramified quadratic extension of $\Q_\ell$, and, for $r>0$ and $i\ge0$, let $U_{\ell^r}^i=\ker(W(\F_{\ell^r})^\times \to W_i(\F_{\ell^r})^\times)$ and let $U_{\ell^r}=W(\F_{\ell^r})^\times$. (Note that $W(\F_{\ell^r})$ here refers to the Witt vectors, not to be confused with the Weyl group.)

Let $\theta_0:\Q_{\ell^2}^\times \to K^\times$ be a character, let $n$ be the smallest integer such that $\theta_0|_{U_{\ell^2}^{n+1}}=1$, and assume that $n$ is positive and even. Assume that $\theta_0|_{U_{\ell^2}^n}$ does not factor through the norm $U_{\ell^2}^n \to U_\ell^n$. This $n$ is the depth of the corresponding supercuspidal representation. Let $\psi:\Z_\ell \to K^\times$ be a character with the property that $n$ is the smallest integer such that $\psi|_{\ell^n\Z_\ell}=1$. There is an element $u \in \Q_{\ell^2}^\times$ such that
\[
\theta_0(1+x) = \psi(\tr_{\Q_{\ell^2}/\Q_\ell}(ux))
\]
for all $x \in \ell^{n/2}W(\F_{\ell^2})$; fix such a $u$.

For an integer $r \ge 0$, let $\Gamma(\ell^n)=\ker\left(\GL_2(\Z_\ell) \to \GL_2(\Z/\ell^r\Z)\right)$. Fix an embedding $U_{\ell^2} \subset \GL_2(\Z_\ell)$, and consider the following subgroups of $\GL_2(\Q_\ell)$:
\begin{align*}
& J''=U_{\ell^2}^1\Gamma(\ell^{n/2+1}) \\
& J'=U_{\ell^2}^1\Gamma(\ell^{n/2}) \\
&J=\Q_{\ell^2}^\times\Gamma(\ell^{n/2})
\end{align*}
Let $\theta:J'' \to K^\times$ be the character defined by
\[
\theta(x(1+\ell y)) = \theta_0(x)\psi(\tr(uy)) 
\]
for all $x \in U_{\ell^2}^1$ and $y \in \ell^{n/2}M_2(\Z_\ell)$; it is well-defined by the definition of $u$. Then $\ker(\theta)$ is a normal subgroup of $J$. 

We can now apply Example \eqref{ex:higher depth} to this situation. Let $H=J'/\ker(\theta)$ and $A=\Q_{\ell^2}^\times/U_{\ell^2}^1\Q_{\ell}^\times \cong \F_{\ell^2}^\times/\F_\ell^\times$, which is a cyclic group of order $\ell+1$, and $G=H \rtimes A$. The center of $H$ is $Z=J''/\ker(\theta)$, on which $\theta$ is a faithful character. The hypotheses of Example \eqref{ex:higher depth} are satisfied by \cite[Section  22.4, pg.~136]{BH2006}. For $R \in \{\cO,K\}$, let $\cC_R(G,\theta_0)$ denote the category of $R[G]$-modules $M$ such that $M|_Z$ is (multiple copies of) $\theta_0$. There is an exact functor
\[
F_R:\cC_R(G,\theta_0) \to R[J/\ker(\theta)]\mathrm{-mod}
\]
defined as follows. Let $\Q_{\ell^2}^\times$ act on $H$ via the quotient map $\Q_{\ell^2}^\times \onto A$, and let $G'=H \rtimes \Q_{\ell^2}^\times$, so that $G$ is a quotient of $G'$. There is a surjective group homomorphism
\[
\phi: G' \to J/\ker(\theta), \ \phi(x,h)=xh
\]
with $\ker(\phi)$ equal to the set of $(x,z) \in  U_{\ell^2}^1 \times Z$ such that $z$ is the image of $x^{-1}$. Let $\pi \in \cC_R(G,\theta_0)$, and let ${\pi'}$ and $\theta'_0$ be the inflations of $\pi$ and $\theta_0$, respectively, to $G'$. Since the central character of $\pi$ is $\theta_0$, the representation $\theta'_0 \otimes_R \pi'$ is trivial on $\ker(\phi)$, so it factors through $\phi$, giving a $R[J/\ker(\theta)]$-module. Define $F_R(\pi)=\theta'_0 \otimes_R \pi'$.

Example \ref{ex:higher depth} defines a representation $\eta_1$ in $\cC_K(G,\theta_0)$; the supercuspidal representation $\pi_{\theta_0}$ associated to $\theta_0$ is defined as the compact induction from $J$ to $\GL_2(\Q_\ell)$ of the inflation of $F_K(\eta_1)$ to $J$. The representation $F_K(\eta_1)$ itself (or sometimes $\Ind_{J \cap \GL_2(\Z_\ell)}^{\GL_2(\Z_\ell)}(F_K(\eta_1)|_{J \cap \GL_2(\Z_\ell)})$) is called the cuspidal type of $\pi_{\theta_0}$. Applying $F_\cO$ to \eqref{eq:eta1} gives a modular analog of the character formula for the cuspidal type of $\pi_{\theta_0}$.
\end{ex}

\section{Minimal lifts}
\label{sec:min lifts}
In this section, we prove Theorem \ref{thm:minimal lifts}. Let $p$, $\ell$, $f$, $\theta$ and $\theta'$ be as in the statement of the theorem. We start by giving a rephrasing of the theorem in terms of cohomology of arithmetic groups. The relationship between the modular forms and the cohomology of arithmetic groups is classical, but may not be familiar to the audience of this paper (especially those interested in modular representation theory). We review some of the ideas, and encourage the interested reader to consult the nice survey in \cite{grossHecke}.

The arithmetic group $\Gamma_0(\ell^2)$ is conjugate to 
\[
\Gamma_0(\ell,\ell) = \left\{ \ttmat{a}{b}{c}{d} \in \SL_2(\Z) \, \middle|\, \ell \mid b,\ \ell \mid c\right\}
\]
which contains $\Gamma(\ell)=\ker(\SL_2(\Z) \to \SL_2(\F_\ell))$ as a normal subgroup. The cohomology of $\Gamma(\ell)$ can be computed using Shapiro's lemma as
\[
H^1(\Gamma(\ell),\C) = H^1(\SL_2(\Z),\C[\SL_2(\F_p)]).
\]
The module $\C[\SL_2(\F_p)]$ decomposes as a direct sum of irreducible representations. For a primitive character $\chi$ of the nonsplit torus, let $\pi_{\chi}$ denote the associated cuspidal representation and $\pi_\chi^\ast$ its $\C$-linear dual. Then there is a modular form of weight 2 and level $\Gamma_0(\ell^2)$ that is supercuspidal at $\ell$ with character $\chi$ if and only if the group $H^1(\SL_2(\Z),\pi_{\chi}^\ast)$ is non-zero. This is discussed in \cite[Sections 9, 10]{grossHecke}; see also \cite[Section 6.1]{CDT}.

In order to consider congruences, we need lattices. Let $\tilde{\pi}_{\chi}$ be an $\cO$-lattice in $\pi_{\chi}$ and $\tilde{\pi}_\chi^\ast$ its $\cO$-linear dual. Let $\bT$ be the $\cO$-subalgebra of $\End_\cO(H^1(\Gamma_0(\ell^2),\cO))$ generated by the Hecke operators $T_n$ for $n\in \Z$ with $3\ell \nmid n$. Let $I_f=\ann_\bT(f)$ be the ideal generated by $T_n-a_n(f)$ for all $n\in \Z$ with $\ell \nmid n$, and let $\m_f$ be the maximal ideal of $\bT$ containing $I_f$. Then Theorem \ref{thm:minimal lifts} is equivalent to the following.
\begin{theorem}
    \label{thm:cohomology min lifts}
    The localization $H^1(\SL_2(\Z),\tilde{\pi}_{\theta'}^\ast)_{\m_f}$ is nonzero.
\end{theorem}

Indeed, if $\bT_{\theta'}$ denotes the quotient of $\bT$ by the annihilator of $H^1(\SL_2(\Z),\tilde{\pi}_{\theta'}^\ast)$, then the minimal prime ideals in $\bT_{\theta'}$ correspond to newforms $g$ that satisfy (2) of Theorem \ref{thm:minimal lifts}. Theorem \ref{thm:cohomology min lifts} implies that there exists such a prime ideal that is contained in $\m_f$, and hence the corresponding form $g$ satisfies (1) of Theorem~\ref{thm:minimal lifts}.

\begin{proof}[Proof of Theorem \ref{thm:cohomology min lifts}]
The long exact sequence in cohomology applied to the short exact sequence
\[
0 \to \tilde{\pi}_{\theta'}^\ast \xrightarrow{\varpi} \tilde{\pi}_{\theta'}^\ast \to \ov \pi_{\theta'}^\ast \to 0
\]
yields an exact sequence
\begin{equation}
\label{eq:theta and theta'}
    0 \to H^1(\SL_2(\Z),\tilde{\pi}_{\theta'}^\ast)_{\m_f} \otimes_\cO k \to H^1(\SL_2(\Z),\ov \pi_{\theta'}^\ast)_{\m_f} \to H^2(\SL_2(\Z),\tilde{\pi}_{\theta'}^\ast)_{\m_f}.
\end{equation}
There is a similar exact sequence for $\tilde{\pi}_\theta^\ast$; note that the assumption $\ov \theta = \ov{\theta'}$ implies that $\ov \pi_{\theta'}^\ast=\ov \pi_{\theta}^\ast$. Since $f$ defines a nonzero class in the localization $H^1(\SL_2(\Z),\tilde{\pi}_{\theta}^\ast)_{\m_f}$, it follows from this similar sequence that $H^1(\SL_2(\Z),\ov \pi_\theta^\ast)_{\m_f}=H^1(\SL_2(\Z),\ov \pi_{\theta'}^\ast)_{\m_f}$ is nonzero. Then, by \eqref{eq:theta and theta'}, it is enough to show that $H^2(\SL_2(\Z),\tilde{\pi}_{\theta'}^\ast)=0$.

Using Theorem \ref{thm:main intro}, we compute $H^2(\SL_2(\Z),\tilde{\pi}_{\chi}^\ast)$ for arbitrary $\chi$. Since there is a finite-index subgroup of $\SL_2(\Z)$ that has cohomological dimension one (indeed, there is such a subgroup that acts freely on the (contractible) upper-half plane), intuitively, the cohomology $\SL_2(\Z)$ in degrees greater than one should all be accounted for by torsion in $\SL_2(\Z)$. 
The theory of Farrell cohomology makes this intuition precise: \cite[Corollary 7.4, pg.~293]{brownbook} implies that there is an isomorphism
\[
H^2(\SL_2(\Z),\tilde{\pi}_\chi^\ast) \cong H^2(P,\tilde{\pi}_\chi^\ast)
\]
where $P \le \SL_2(\Z)$ is a subgroup of order $3$ (which is unique up to conjugacy). Here, $P$ acts on $\tilde{\pi}_\chi$ 
through its image in $G=\GL_2(\F_\ell)$, which is a subgroup of a non-split maximal torus $T$ in $G$; we denote the image of $P$ in $T$ by $P$ (abusing notation somewhat).

By Theorem \ref{thm:main intro} (see also \eqref{eq:GL2 sequence}), there is an exact sequence of $\cO[G]$-modules
\[
0 \to \tilde{\pi}_\chi \to e\Ind_{ZU}^G \tilde{\psi} \to e \Ind_T^G \tilde{\chi} \to 0
\]
where $e\Ind_{ZU}^G \tilde{\psi}$ is projective (and therefore cohomologically trivial). 
Taking $\cO$-linear duals (which preserves projectivity and, because the modules are $\cO$-free, exactness), there is an exact sequence
\[
0 \to  (e \Ind_T^G \tilde{\chi})^* \to (e\Ind_{ZU}^G \tilde{\psi})^* \to \tilde{\pi}_\chi^*\to 0
\]
where the middle term is cohomologically trivial. There are thus isomorphisms
\begin{equation}
\label{eq:H^2 H^3 shift}
    H^2(\SL_2(\Z),\tilde{\pi}_\chi^\ast) \cong H^2(P,\tilde{\pi}_\chi^\ast) \cong H^3(P,(e\Ind_T^G\tilde\chi)^\ast).
\end{equation}

To compute $H^3(P,(e\Ind_T^G\tilde\chi)^\ast)$, we need to analyze the restriction of $e\Ind_T^G\tilde\chi$ to $P$. Mackey's formula implies that this restriction is
\[
\left(\Ind_T^G\tilde\chi\right)|_P=\bigoplus_{\bar s \in T \backslash G /P} \Ind_{s^{-1}Ts \cap {P}}^P \tilde{\chi}^s,
\]
where $s \in G$ is an element of the double coset $\bar s$ and where $\tilde{\chi}^s$ is the conjugate character $\tilde{\chi}^s(x)=\tilde{\chi}(sxs^{-1})$.
Since $P$ is generated by a regular element of $T$, the subgroup $s^{-1}Ts \cap {P}$ is equal to $P$ if and only if $\bar s$ is the double coset of an element of the Weyl group of $T$ (and is the trivial group otherwise). It follows that the restriction of $\Ind_T^G\tilde\chi$ to $P$ is $\tilde \chi|_P \oplus \tilde \chi^w|_P \oplus F$ for a free $\cO[P]$ module $F$. Since $e\Ind_T^G\tilde\chi$ is a direct summand of $\Ind_T^G\tilde\chi$, it follows from this and a comparison of trace characters that 
\begin{equation}
    \label{eq: resP of Ind}
    \left(e\Ind_T^G\tilde\chi\right)|_P =\tilde \chi|_P \oplus \tilde \chi^w|_P \oplus M
\end{equation}
for a projective $\cO[P]$-module $M$.

Together, \eqref{eq:H^2 H^3 shift} and \eqref{eq: resP of Ind} show that there is an isomorphism
\[
H^2(\SL_2(\Z),\tilde{\pi}_\chi^\ast) \cong H^3(P,\tilde\chi|_P^{-1} \oplus \tilde\chi^w|_P^{-1}).
\]
By $2$-periodicity of cohomology of cyclic groups (see \cite[Section VI.9, pg.~154]{brownbook}), for a character $\lambda:P \to \cO^\times$, the group $H^3(P,\lambda) \cong \hat H^{-1}(P,\lambda)$ is zero if and only if $\lambda$ is trivial. Hence $H^2(\SL_2(\Z),\tilde{\pi}_\chi^\ast)$ is zero if and only if $\tilde\chi$ has prime-to-$3$ order. Since $\theta'$ has prime-to-$3$ order by definition, this proves $H^2(\SL_2(\Z),\tilde{\pi}_{\theta'}^\ast)=0$.
\end{proof}

To see how this proof breaks down in the setting of Example \ref{ex:3 11}, note that the above computation shows that $H^2(\SL_2(\Z),\tilde{\pi}_\chi^\ast)$ is nonzero if $\chi|_P \ne 1$. This means that, without the assumption that the order of $\theta'$ is prime-to-$3$, the group $H^2(\SL_2(\Z),\tilde{\pi}_{\theta'}^\ast)_{\m_f}$ could be nonzero and provide an obstruction to lifting in \eqref{eq:theta and theta'}. This is what happens in Example \ref{ex:3 11}.

%%%%%%%%%%%%%%%%%%%%%%%%%%%% References %%%%%%%%%%%%%%%%%%%%%%%%%%%%%%

\bibliography{bibliography1}{}
\bibliographystyle{alpha}

\end{document}